\documentclass[11pt,reqno]{amsart}
\usepackage{amsthm}
\usepackage{amssymb}
\usepackage{graphics}
\usepackage{tikz}
\usetikzlibrary{shapes,backgrounds,calc}
\usepackage{latexsym}
\usepackage{multicol}
\usepackage{verbatim,enumerate}
\usepackage{accents}
\usepackage{cite}

\usepackage[colorlinks=true, linkcolor=blue, citecolor=blue, urlcolor=blue]{hyperref}
\usepackage{hyperref}
\usepackage{amsmath, amscd,url}

\advance\textwidth by 1.3in \advance\oddsidemargin by -.6in \advance\evensidemargin by -.6in
\theoremstyle{definition}
\newtheorem{lem}{Lemma}
\newtheorem{prop}{Proposition}

\theoremstyle{definition}
\newtheorem{defn}{Definition}

\theoremstyle{definition}

\newcounter{cnt}
 \makeatletter
\def\mydggeometry{\makeatletter\dg@YGRID=1\dg@XGRID=20\unitlength=0.003pt\makeatother}
\makeatother \theoremstyle{remark}

\numberwithin{equation}{section}
\let\bwdg\bigwedge
\def\bigwedge{{\textstyle\bwdg}}

\newcommand{\Evec}{\vec{E}}

\newcommand{\nc}{\newcommand}
\newcommand{\rnc}{\renewcommand}

\nc{\cal}{\mathcal} \nc{\goth}{\mathfrak} \rnc{\bold}{\mathbf}

\nc\bomega{{\mbox{\boldmath $\omega$}}} \nc\bpsi{{\mbox{\boldmath $\Psi$}}}
 \nc\balpha{{\mbox{\boldmath $\alpha$}}}
 \nc\bpi{{\mbox{\boldmath $\pi$}}}
 \nc\bvpi{{\mbox{\boldmath $\varpi$}}}
\nc\chara{\operatorname{ch}}

  \nc\bxi{{\mbox{\boldmath $\xi$}}}
\nc\bmu{{\mbox{\boldmath $\mu$}}} \nc\bcN{{\mbox{\boldmath $\cal{N}$}}} \nc\bcm{{\mbox{\boldmath $\cal{M}$}}} \nc\blambda{{\mbox{\boldmath
$\lambda$}}}\nc\bnu{{\mbox{\boldmath $\nu$}}}

\makeatletter
\def\section{\def\@secnumfont{\mdseries}\@startsection{section}{1}%
  \z@{.7\linespacing\@plus\linespacing}{.5\linespacing}%
  {\normalfont\scshape\centering}}
\def\subsection{\def\@secnumfont{\bfseries}\@startsection{subsection}{2}%
  {\parindent}{.5\linespacing\@plus.7\linespacing}{-.5em}%
  {\normalfont\bfseries}}
\makeatother

 \nc{\Hom}{\operatorname{Hom}}
  \nc{\mode}{\operatorname{mod}}
\nc{\End}{\operatorname{End}} \nc{\wh}[1]{\widehat{#1}} \nc{\Ext}{\operatorname{Ext}} \nc{\ch}{\text{ch}} \nc{\ev}{\operatorname{ev}}
\nc{\Ob}{\operatorname{Ob}} \nc{\soc}{\operatorname{soc}} \nc{\rad}{\operatorname{rad}} \nc{\head}{\operatorname{head}}

 \nc{\Cal}{\cal} \nc{\Xp}[1]{X^+(#1)} \nc{\Xm}[1]{X^-(#1)}
\nc{\on}{\operatorname} \nc{\Z}{{\bold Z}} \nc{\J}{{\cal J}}  \nc{\Q}{{\bold Q}}

\nc{\N}{{\bold N}}  \nc\boa{\bold a} \nc\bob{\bold b} \nc\boc{\bold c} \nc\bod{\bold d} \nc\boe{\bold e} \nc\bof{\bold f} \nc\bog{\bold g}
\nc\boh{\bold h} \nc\boi{\bold i} \nc\boj{\bold j} \nc\bok{\bold k} \nc\bol{\bold l} \nc\bom{\bold m} \nc\bon{\mathbb n} \nc\boo{\bold o}
\nc\bop{\bold p} \nc\boq{\bold q} \nc\bor{\bold r} \nc\bos{\bold s} \nc\boT{\bold t} \nc\boF{\bold F} \nc\bou{\bold u} \nc\bov{\bold v}
\nc\bow{\bold w} \nc\boz{\bold z}\nc\ba{\bold A} \nc\bb{\bold B} \nc\bc{\mathbb C} \nc\bd{\bold D} \nc\be{\bold E} \nc\bg{\bold
G} \nc\bh{\bold H} \nc\bi{\bold I} \nc\bj{\bold J} \nc\bk{\bold K} \nc\bl{\bold L} \nc\bm{\bold M} \nc\bn{\mathbb N} \nc\bo{\bold O} \nc\bp{\bold
P} \nc\bq{\bold Q} \nc\br{\bold R} \nc\bs{\bold S} \nc\bt{\bold T} \nc\bu{\bold U} \nc\bv{\bold V} \nc\bw{\bold W} \nc\bz{\mathbb Z} \nc\bx{\bold
x} \nc\KR{\bold{KR}} \nc\rk{\bold{rk}} \nc\het{\text{ht }}

\nc\toa{\tilde a} \nc\tob{\tilde b} \nc\toc{\tilde c} \nc\tod{\tilde d} \nc\toe{\tilde e} \nc\tof{\tilde f} \nc\tog{\tilde g} \nc\toh{\tilde h}
\nc\toi{\tilde i} \nc\toj{\tilde j} \nc\tok{\tilde k} \nc\tol{\tilde l} \nc\tom{\tilde m} \nc\ton{\tilde n} \nc\too{\tilde o} \nc\toq{\tilde q}
\nc\tor{\tilde r} \nc\tos{\tilde s} \nc\toT{\tilde t} \nc\tou{\tilde u} \nc\tov{\tilde v} \nc\tow{\tilde w} \nc\toz{\tilde z} \nc\woi{w_{\omega_i}}

\begin{document}


\title{The Peterson recurrence formula for the chromatic discriminant of a graph }
\author{G. Arunkumar}
\address{The Institute of Mathematical Sciences, HBNI, Chennai, India}
\email{gakumar@imsc.res.in}
\thanks{The author acknowledges partial support under the DST Swarnajayanti fellowship DST/SJF/MSA-02/2014-15 of Amritanshu Prasad.}


\begin{abstract}
	The absolute value of the coefficient of $q$ in the chromatic polynomial of a graph $G$ is known as the {\em chromatic discriminant} of $G$ and is denoted $\alpha(G)$. There is a well known recurrence formula for $\alpha(G)$ that comes from the deletion-contraction rule for the chromatic polynomial. In this paper we prove another recurrence formula for $\alpha(G)$ that comes from the theory of Kac-Moody Lie algebras. We start with a brief survey on many interesting algebraic and combinatorial interpretations of $\alpha(G)$. We use two of these interpretations (in terms of  acyclic orientations and spanning trees) to give two bijective proofs for our recurrence formula of $\alpha(G)$. 
\end{abstract}

\maketitle

\section{Introduction}

Let $G$ be a simple graph and let $\chi(G,q)$ denote its chromatic polynomial. The absolute value of the coefficient of $q$ in $\chi(G,q)$ is known as the {\em chromatic discriminant} of the graph $G$ \cite{MR1861053,shi2016graph} and is denoted  $\alpha(G)$.  It is an important graph invariant with numerous algebraic and combinatorial interpretations.
For instance, letting $q$ denote a fixed vertex of the graph $G$, it is well known that each of the following sets has cardinality $\alpha(G)$:
\begin{enumerate}
\item Acyclic orientations of $G$ with unique sink at $q$ \cite{GZ83},
\item Maximum $G$-parking functions relative to $q$ \cite{MR2592488},
\item Minimal $q$-critical states \cite[Lemmas 14.12.1 and 14.12.2]{MR1829620},
\item Spanning trees of $G$ without broken circuits \cite{MR830053},
\item Conjugacy classes of Coxeter elements in the Coxeter group associated to $G$ \cite{MR2515767,MR1857933,MR2431028},
\item Multilinear Lyndon heaps on $G$ \cite{lalonde,MR1235180, viennot-imsc}.
  \end{enumerate}

 In addition, $\alpha(G)$ is also equal to the dimension of the root space corresponding to the sum of all simple roots in the Kac-Moody Lie algebra  associated to $G$  \cite{MR3342717,akv}.

 We have the following recurrence formula for $\alpha(G)$ (see for instance \cite{MR1778205}) which is an immediate consequence of the well-known {\em deletion-contraction rule} for the chromatic polynomial: 
 \begin{equation}\label{eq:dcrecc}
 	\alpha(G) = \alpha(G\backslash e) + \alpha(G/e),
 \end{equation}
  where $e$ is any edge of $G$. Here, $G\backslash e$ denotes $G$ with $e$ deleted and $G/e$ denotes the simple graph obtained from $G$ by identifying the two ends of $e$ ({\em i.e., contracting $e$ to a single vertex}) and removing any multiple edges that result.

 Yet another recurrence formula for $\alpha(G)$ was obtained in \cite{MR3342717} using its connection to root multiplicities of Kac-Moody Lie algebras. To state this, we introduce some notation: for a graph $G$, we let $V(G)$ and $E(G)$ denote its vertex and edge sets respectively. We say that the ordered pair $(G_1, G_2)$ is an {\em ordered partition of  $G$}, if $G_1$ and $G_2$ are non-empty subgraphs of $G$ whose vertex sets form a partition of $V(G)$, i.e., they are disjoint and their union is $V(G)$. When we don't care about the ordering of $G_1, G_2$, we call the set $\{G_1,G_2\}$ an {\em unordered partition of $G$}. We say that an edge $e$ {\em straddles} $G_1$ and $G_2$ if one end of $e$ is in $G_1$ and the other in $G_2$.
 
 We then have:
 \begin{prop} \cite{MR3342717}
 \begin{equation}\label{eq:petrecc}
 	 2~e(G)~ \alpha(G) = \sum\limits_{\substack{(G_1,G_2) \\ \text{ordered partitions of } G}} \alpha(G_1)~\alpha(G_2)~e(G_1,G_2).
 \end{equation}
 Here $e(G)$ is the total number of edges in $G$, $e(G_1,G_2)$ is the number of edges that straddle $G_1$ and $G_2$, and the sum ranges
 over ordered partitions of $G$.
\end{prop}

 We note that the recurrence formula \eqref{eq:petrecc} does not seem to follow directly from \eqref{eq:dcrecc}.
 In \cite{MR3342717}, \eqref{eq:petrecc} was derived from the {\em Peterson recurrence formula} \cite{kac} for root multiplicities of Kac-Moody Lie algebras. The goal of this paper is to give a purely combinatorial (bijective) proof of \eqref{eq:petrecc}.
 
 To construct a bijective proof, we need sets whose cardinalities are the left and right hand sides of \eqref{eq:petrecc}. We in fact give two bijective proofs, starting
 with  the interpretations of $\alpha(G)$ in terms of acyclic orientations and spanning trees.
 
 \noindent
 {\em Acknowledgements}$\colon$ The author would like to thank Sankaran Viswanath for many fruitful discussions.


\section{Acyclic orientations with unique fixed sink}

In this section we give a bijective proof of the recurrence formula \eqref{eq:petrecc}  in terms of acyclic orientations. 


We recall that an acyclic orientation of $G$ is an assignment of arrows to its edges such that there are no directed cycles in the resulting directed graph.
A sink in an acyclic orientation is a vertex which only has incoming arrows. The set of all acyclic orientations of $G$ is denoted $\mathcal{A}(G)$. For a vertex $q$ of
$G$, the set of all acyclic orientations in which $q$ is the unique sink is denoted $\mathcal{A}(G,q)$. It is well-known that the cardinality of $\mathcal{A}(G,q)$ is independent of $q$ and equals $\alpha(G)$ \cite{GZ83}. 
The following characterization of $\mathcal{A}(G,q)$ is immediate.

\begin{lem} \label{Lemma:uniq_sink_char}
	Fix a vertex $q$ of $G$ and let $\lambda \in \mathcal{A}(G)$. Then $\lambda \in \mathcal{A}(G,q)$ if and only if for every $p \in V(G)$, there is a directed path in $\lambda$ from $p$ to $q$.
\end{lem}

This motivates the following:

\begin{defn} \label{def:reachables}
  Given a vertex $q$ and an acyclic orientation $\lambda$ of $G$, let $V(\lambda,q)$ denote the set of all vertices $p$ for which there is a directed path in $\lambda$ from $p$ to $q$. We call this the set of $q$-reachable vertices in $\lambda$.
\end{defn}

  We record the following simple observation:

\begin{lem} \label{Lemma:nbd_char}
  Let $x$ be a vertex of $G$.
  \begin{enumerate}
    \item[(a)] If $x \not\in V(\lambda,q)$, then $x \not\in V(\lambda,p)$ for all $p \in V(\lambda,q)$. \\
    \item[(b)] In particular, an edge joining $p$ and $x$ with $p \in V(\lambda,q)$ and $x \not\in V(\lambda,q)$ is directed from $p$ to $x$ in $\lambda$.
  \end{enumerate}
\end{lem}

Our next goal is to construct sets $A$ and $B$ whose cardinalities are respectively equal to the left and right hand sides of \eqref{eq:petrecc} and to exhibit a bijection between them. To this end, we first consider the set $\Evec$ of oriented edges of $G$; an element of $\Evec$ is an edge
of $G$ with an arrow marked on it (in one of two possible ways). Thus $\Evec$ has cardinality $2\,e(G)$. If $\vec{e}$ is an element of $\Evec$ corresponding to an edge joining vertices $p$ and $q$ with the arrow pointing from $p$ to $q$, we call $p$ the {\em tail} of $\vec{e}$ and $q$ its {\em head}.

We now define $A$ to be the set consisting of pairs $(\vec{e}, \lambda) \in \Evec \times \mathcal{A}(G)$ such that the head of $\vec{e}$ is the unique sink of $\lambda$. For a fixed $\vec{e}$, there are $\alpha(G)$ choices for $\lambda$ since $\lambda$ ranges over $\mathcal{A}(G,q)$ where $q$ is the head of $\vec{e}$. It is now clear that $A$ has cardinality exactly $2\,e(G)\,\alpha(G)$.

To define $B$, we first take an ordered partition $(G_1, G_2)$ of $G$. Let $E(G_1,G_2)$ denote the set of edges
straddling $G_1$ and $G_2$. Let $B(G_1,G_2)$ denote the set of triples $(e,\lambda_1, \lambda_2)$ where $e \in E(G_1,G_2)$, say $e$ joins $p_1$ and $p_2$ with $p_i$ a vertex of $G_i$, $i=1,2$, and $\lambda_i$ is an acyclic orientation of $G_i$ with unique sink at $p_i$, $i=1,2$. Arguing as before, one concludes that $B(G_1,G_2)$ has cardinality $\alpha(G_1)\, \alpha(G_2)\, e(G_1,G_2)$. We now let $B$ denote the disjoint union of the $B(G_1, G_2)$ over all ordered partitions $(G_1,G_2)$ of $G$. It clearly has cardinality equal to the right hand side of \eqref{eq:petrecc}.

We now define a map $\varphi: A \to B$ which will turn out to be the bijection we seek. Given $(\vec{e}, \lambda) \in A$, let $p$ and $q$ denote the tail and head of $\vec{e}$ respectively. Note that $\lambda \in \mathcal{A}(G,q)$. Let $V_1 = V(\lambda,p)$ denote the set of $p$-reachable vertices in $\lambda$ (definition \ref{def:reachables}) and let $V_2 = V(G) \backslash V_1$. Observe that $p \in V_1$ and $q \in V_2$. For $i=1,2$, let $G_i$ denote the subgraphs of $G$ induced by $V_i$, and let $\lambda_i$ denote the restriction of $\lambda$ to $G_i$.

We claim that $\lambda_1$ has a unique sink at $p$ and $\lambda_2$ has a unique sink at $q$. The first assertion follows simply from Lemma \ref{Lemma:uniq_sink_char}. For the second assertion, observe that if $x \in V_2 \subset V(G)$, then there is a directed path from $x$ to $q$ in $\lambda$. Since $x \notin V(\lambda,p)$, Lemma \ref{Lemma:nbd_char}(a) implies that no vertex of this directed path can lie in $V_1$. In other words this directed path is entirely within $G_2$, and we are again done by Lemma \ref{Lemma:uniq_sink_char}.

Let $e$ denote the undirected edge joining $p$ and $q$. We have thus shown that the triple $(e,\lambda_1,\lambda_2)$ is in $B(G_1,G_2) \subset B$.
We define $\varphi(\vec{e}, \lambda) =  (e,\lambda_1,\lambda_2)$.

To see that $\varphi$ is a bijection, we describe its inverse map. Let $(G_1,G_2)$ be an ordered partition of $G$; given a triple $(e,\lambda_1,\lambda_2) \in B(G_1,G_2)$, we construct an acyclic orientation $\lambda$ of $G$ as follows: on $G_1$ and $G_2$, we define $\lambda$ to coincide with $\lambda_1$ and $\lambda_2$ respectively. It only remains to define an orientation for the straddling edges (this includes $e$); we orient all of them pointing from $G_1$ towards $G_2$, i.e., such that their tails are in $G_1$ and their heads in $G_2$. We let $\vec{e}$ denote the edge $e$ with the above orientation.

We claim $(\vec{e}, \lambda) \in A$. First observe that $\lambda$ is in fact acyclic; since $\lambda$ extends $\lambda_i$ for $i=1,2$, any directed cycle of $\lambda$ must necessarily involve vertices from both $G_1$ and $G_2$. But this is impossible since all straddling edges point the same way, from $G_1$ towards $G_2$. Let $p, q$ denote the tail and head of $\vec{e}$. It remains to show that $\lambda$ has a unique sink at $q$, or equivalently, by Lemma \ref{Lemma:uniq_sink_char}, that there is a directed path in $\lambda$ from any vertex $x$ to $q$. This is clear if $x$ is a vertex of $G_2$. If $x$ is in $G_1$, we have a directed path in $\lambda_1$ from $x$ to $p$. Now the edge $\vec{e}$ is directed from $p$ to $q$; concatenating this directed path with $\vec{e}$ produces a directed path from $x$ to $q$ in $\lambda$ as required. We define the map $\psi: B \to A$ by $\psi(e,\lambda_1,\lambda_2) = (\vec{e}, \lambda)$.

Observe that for the $\lambda$ defined above, the set of $p$-reachable vertices is exactly $V(G_1)$. This is because edges straddling $G_1$ and $G_2$ point away from $G_1$, so no vertex of $G_2$ is $p$-reachable. This implies that $\varphi \circ \psi$ is the identity map on $B$. Further, it readily follows from Lemma \ref{Lemma:nbd_char}(b) that $\psi \circ \varphi$ is the identity map on $A$. This establishes that $\varphi$ is a bijection. \qed

\section{Spanning trees without broken circuits}

In this section we give another bijective proof of the recurrence formula \eqref{eq:petrecc}, this time using the fact that $\alpha(G)$ counts the number of spanning trees of $G$ without broken circuits.


\begin{defn}
Let $\sigma$ be a total ordering on the set $E(G)$ of edges of $G$. Given a circuit in $G$, it has a unique maximum edge with respect to $\sigma$; the set of edges obtained by deleting this edge from the circuit is called a {\em broken circuit} relative to $\sigma$. The set of all broken circuits relative to $\sigma$ is denoted $B_G(\sigma)$.
\end{defn}

Let $S_G(\sigma)$ be the set of all spanning trees of $G$ that contain no broken circuits relative to $\sigma$. It is well-known that the cardinality of $S_G(\sigma)$ is independent of the choice of $\sigma$, and equals $\alpha(G)$ \cite{MR830053}.

Given a total ordering $\sigma$ on $E(G)$, let $\max(\sigma)$ denote the maximum edge in $E(G)$.
The following lemma is immediate.

\begin{lem} \label{Lemma:sptreemaxedge}
  Any spanning tree in $S_G(\sigma)$ contains the edge $\max(\sigma)$.
\end{lem}

In the sequel, we will fix for each edge $e$,  a total order $\sigma_e$ on $E(G)$ for which $\max(\sigma_e) = e$. We will write $B_G(e)$ and $S_G(e)$ for the sets $B_G(\sigma_e)$ and $S_G(\sigma_e)$ respectively.

We now proceed to prove \eqref{eq:petrecc} in the following equivalent form:
\begin{equation} \label{eq:petreccprime}
e(G)\, \alpha(G) = \sum\limits_{\substack{\{G_1,G_2\} \\ \text{unordered partitions of } G}} \alpha(G_1) \, \alpha(G_2) \, e(G_1,G_2).
\end{equation}

We first define the set $A$ to consist of pairs $(e, T)$ where $e$ is an edge and $T \in S_G(e)$; from the above discussion, $A$ has cardinality $e(G)\,\alpha(G)$.

Next, we define the set $B$. Given an unordered partition $\{G_1,G_2\}$ of $G$, define $B(\{G_1,G_2\})$ to be the set of pairs $(e,\{T_1,T_2\})$ where $e$ is an edge that straddles $G_1$ and $G_2$ and $T_i \in S_{G_i}(e)$ for $i=1,2$. Here, $S_{G_i}(e)$ is the set of spanning trees of $G_i$ which contain no broken circuits relative to the total order $\sigma_e$ restricted to the edges of $G_i$.
We let $B$ denote the disjoint union of $B(\{G_1,G_2\})$ as $\{G_1,G_2\}$ ranges over unordered partitions of $G$. Clearly $B$ has cardinality equal to the right hand side of \eqref{eq:petreccprime}.

We define maps $\varphi: A \to B$ and $\psi: B \to A$ as follows:

Given $(e, T) \in A$, $e$ occurs in $T$ in view of Lemma \ref{Lemma:sptreemaxedge}. Deleting $e$ from the spanning tree $T$ will result in a pair of trees $T_1, T_2$ with vertex sets $V_1$ and $V_2$. Let $G_i$ denote the subgraph induced by $V_i$, $i=1,2$; clearly $\{G_1,G_2\}$ is an unordered partition of $G$ and $e$ straddles the $G_i$. Observe that since the total order on $E(G_i)$ is defined as the restriction of the total order $\sigma_e$ on $E(G)$, $T_i$ will contain no broken circuits of $G_i$ for $i=1,2$, i.e., $T_i \in S_{G_i}(e)$. We set $\varphi(e,T) = (e,\{T_1, T_2\})$.

For the inverse map $\psi$, let $(e,\{T_1, T_2\}) \in B$. Define $T$ to be the spanning tree of $G$ obtained by adding the edge $e$ to the union of $T_1$ and $T_2$. To prove that $T$ contains no broken circuits relative to $\sigma_e$, observe that any broken circuit of $T$ cannot lie entirely within $T_1$ or $T_2$, and must hence contain the edge $e$. But $e$ is the maximum edge relative to $\sigma_e$, so this cannot be a broken circuit by definition. Thus $(e,T) \in A$, and we define $\psi(e,\{T_1, T_2\}) = (e,T)$.

It is straightforward to check that $\varphi$ and $\psi$ are indeed inverse maps. \qed.

\bibliographystyle{plain}
\bibliography{Arunkumar}

\def\Dbar{\leavevmode\lower.6ex\hbox to 0pt{\hskip-.23ex \accent"16\hss}D}
\begin{thebibliography}{10}

\bibitem{akv}
G.~{Arunkumar}, D.~{Kus}, and R.~{Venkatesh}.
\newblock {Root multiplicities for Borcherds algebras and graph coloring}.
\newblock {\em arXiv:1612.01320}, December 2016.

\bibitem{MR2592488}
Brian Benson, Deeparnab Chakrabarty, and Prasad Tetali.
\newblock {$G$}-parking functions, acyclic orientations and spanning trees.
\newblock {\em Discrete Math.}, 310(8):1340--1353, 2010.

\bibitem{MR830053}
Andreas Blass and Bruce~Eli Sagan.
\newblock Bijective proofs of two broken circuit theorems.
\newblock {\em J. Graph Theory}, 10(1):15--21, 1986.

\bibitem{MR2515767}
Henrik Eriksson and Kimmo Eriksson.
\newblock Conjugacy of {C}oxeter elements.
\newblock {\em Electron. J. Combin.}, 16(2, Special volume in honor of Anders
  Bj\"orner):Research Paper 4, 7, 2009.

\bibitem{MR1778205}
David~D. Gebhard and Bruce~E. Sagan.
\newblock Sinks in acyclic orientations of graphs.
\newblock {\em J. Combin. Theory Ser. B}, 80(1):130--146, 2000.

\bibitem{MR1829620}
Chris Godsil and Gordon Royle.
\newblock {\em Algebraic graph theory}, volume 207 of {\em Graduate Texts in
  Mathematics}.
\newblock Springer-Verlag, New York, 2001.

\bibitem{GZ83}
Curtis Greene and Thomas Zaslavsky.
\newblock On the interpretation of {W}hitney numbers through arrangements of
  hyperplanes, zonotopes, non-{R}adon partitions, and orientations of graphs.
\newblock {\em Trans. Amer. Math. Soc.}, 280(1):97--126, 1983.

\bibitem{kac}
V.~G. Kac.
\newblock {\em Infinite dimensional {L}ie algebras}.
\newblock Cambridge University Press, third edition, 1990.

\bibitem{MR1235180}
Pierre Lalonde.
\newblock Bases de {L}yndon des alg\`ebres de {L}ie libres partiellement
  commutatives.
\newblock {\em Theoret. Comput. Sci.}, 117(1-2):217--226, 1993.
\newblock Conference on Formal Power Series and Algebraic Combinatorics
  (Bordeaux, 1991).

\bibitem{lalonde}
Pierre Lalonde.
\newblock Lyndon heaps: an analogue of {L}yndon words in free partially
  commutative monoids.
\newblock {\em Discrete Math.}, 145(1-3):171--189, 1995.

\bibitem{MR1861053}
Bodo Lass.
\newblock Orientations acycliques et le polyn\^ome chromatique.
\newblock {\em European J. Combin.}, 22(8):1101--1123, 2001.

\bibitem{MR2431028}
Matthew Macauley and Henning~S. Mortveit.
\newblock On enumeration of conjugacy classes of {C}oxeter elements.
\newblock {\em Proc. Amer. Math. Soc.}, 136(12):4157--4165, 2008.

\bibitem{MR1857933}
Jian-yi Shi.
\newblock Conjugacy relation on {C}oxeter elements.
\newblock {\em Adv. Math.}, 161(1):1--19, 2001.

\bibitem{shi2016graph}
Y.~Shi, M.~Dehmer, X.~Li, and I.~Gutman.
\newblock {\em Graph Polynomials}.
\newblock Discrete Mathematics and its applications. CRC Press, 2016.

\bibitem{MR3342717}
R.~Venkatesh and Sankaran Viswanath.
\newblock Chromatic polynomials of graphs from {K}ac-{M}oody algebras.
\newblock {\em J. Algebraic Combin.}, 41(4):1133--1142, 2015.

\bibitem{viennot-imsc}
G\'erard~Xavier Viennot.
\newblock Commutations and heaps of pieces, chapter 5.
\newblock Lectures at IMSc, Chennai.
  \url{http://www.xavierviennot.org/coursIMSc2017/Ch_5_files/cours_IMSc17_Ch5a.pdf}.

\end{thebibliography}

\end{document}